\documentclass[11pt,letterpaper]{amsart}

\usepackage{tikz}
\usetikzlibrary{arrows.meta}
\usepackage{url}
\usepackage{graphicx}

\usepackage{amssymb}
\usepackage[T1]{fontenc}
\usepackage{mathscinet}
\usepackage{enumitem}

\allowdisplaybreaks%

{\theoremstyle{plain}%
 \newtheorem{theorem}{Theorem}
 
 \newtheorem{lemma}{Lemma}%
}
{\theoremstyle{remark}

}
{\theoremstyle{definition}

\newtheorem{example}{Example}
}

\title[The Erd\H{o}s--S\'{a}rk\"{o}zy product-density problem]{A solution to the Erd\H{o}s--S\'{a}rk\"{o}zy product-density problem}
\author{John M.\ Campbell}
\address{Department of Mathematics and Statistics, Dalhousie University,
Halifax, NS B3H 4R2, Canada}
\email{jh241966@dal.ca}
 


\begin{document}

\begin{abstract}
 Color each positive integer red or blue, and call an integer a \emph{same-color product} if it is the product of two red numbers or of two 
 blue numbers. Does there exist a universal constant $c > 0$ such that: For each fixed coloring, and for all sufficiently large $M$ (where 
 ``sufficiently large'' may depend on the coloring), at least $cM$ of the integers among $1$, $2$, $\ldots$, $M$ are same-colored 
 products (allowing for equal factors, and where each distinct product is counted once)? This problem was given (in an equivalent way) in 
 Guy's famous text on unsolved problems in number theory, and was posed (in equivalent ways) by Erd\H{o}s and S\'{a}rk\"{o}zy in 1990, 
 in the proceedings of the 1988 Banff number theory conference, and 1993 [\emph{Colloq.\ Math.}], and this problem seems to have 
 remained open. We succeed in solving this problem, by proving that no such universal constant exists. 
 \end{abstract}

\maketitle

\section{Introduction}\label{secintro}
 Multiplicative Ramsey theory is a deep and challenging area in combinatorial number theory. A foundational reference in multiplicative 
 Ramsey theory is due to Hindman \cite{Hindman1979}, who introduced and proved a result that gives that each finite coloring of 
 $\mathbb N$ admits a strictly increasing infinite sequence such that each nonempty finite product of distinct terms has the same color. 
 Some similarly notable works are due to Bergelson et al.\ \cite{BeiglbockBergelsonHindmanStrauss2008,Bergelson2005}, Green and 
 Sanders \cite{GreenSanders2016}, and Moreira \cite{Moreira2017}. We solve, building on our extensive interactions with GPT-6 Astra Pro, 
 a problem in the field of multiplicative Ramsey theory that is given in Guy's classic textbook on unsolved problems in number theory 
 \cite[p.\ 204]{Guy1994} and that was repeatedly posed by Erd\H os and S\'ark\"ozy \cite{ErdosSarkozy1990,ErdosSarkozy1993}. It appears 
 that this problem has, prior to our work, remained open. 

 Erd\H os and S\'ark\"ozy \cite{ErdosSarkozy1990,ErdosSarkozy1993}
 pose the following question: 
 For each \(\varepsilon>0\), is it true that there is a partition \(\{\mathcal{A}_1,\mathcal{A}_2\}\) of \(\mathbb{N}\) 
 such that the lower asymptotic density of the sequence consisting 
 of the distinct monochromatic products \(aa'\) is less than \(\varepsilon\)?
 In this context, a positive integer $n$ is said to be a \emph{monochromatic product} if 
 it admits at least one factorization of the form $n = a a'$ such that the two factors $a$ and $a'$ 
 are either both in $\mathcal{A}_{1}$ or are both in $\mathcal{A}_{2}$. 
 In this paper, we answer 
 Erd\H os and S\'ark\"ozy's question in the affirmative. 
 The problem given by determining whether or not this question holds
 is equivalent to a problem in Guy's classic text on unsolved problems in number theory \cite[p.\ 204]{Guy1994}, 
 in the sense that a solution to one problem gives a solution to the other problem. More specifically, 
 a negative solution to Guy's problem is equivalent to a solution in the affirmative to the 
 Erd\H os--S\'ark\"ozy problem. 
 For every fixed finite partition, the lower density of the monochromatic products is positive
 \cite{ErdosSarkozy1990}, but 
 what is at issue is whether that positive lower bound can be chosen independently of the partition. 

 For a $2$-coloring $\mathcal{C}$ of $\mathbb{N}$, write
\begin{equation}\label{defineRCN} 
 r_{\mathcal{C}}(N) = \frac{\text{$\#$ distinct monochromatic products among $1$, $2$, $\ldots$, $N$}}{N} 
\end{equation} 
 and 
\begin{equation}\label{definedcalC} 
 d_{\mathcal{C}} = \liminf_{N \to \infty} r_{\mathcal{C}}(N). 
\end{equation}
 Thus the Erd\H os--S\'ark\"ozy problem above is given by whether or not it is the case that 
\begin{equation}\label{ESformulation} 
 \forall \varepsilon > 0 \, \exists \mathcal{C} \ d_{\mathcal{C}} < \varepsilon. 
\end{equation}
 As suggested in the above abstract, the problem given by whether or not \eqref{ESformulation} 
 holds is equivalent to the problem as to whether or not the negation
\begin{equation}\label{abstractformulation} 
 \exists \varepsilon > 0 \, \forall \mathcal{C} \ d_{\mathcal{C}} \geq \varepsilon 
\end{equation}
 holds. 
 More precisely, and in view of \eqref{definedcalC}, we find that 
 \eqref{abstractformulation} is equivalent (as clarified below) to 
\begin{multline*}
 \exists c > 0 \, \forall \mathcal{C} \, \exists N_{\mathcal{C}} \, \forall M \geq N_{\mathcal{C}} \ 
 \text{$\#$\big(distinct monochromatic} \\ 
 \text{products among $1$, $2$, $\ldots$, $M$\big)} 
 \geq c M, 
\end{multline*}
 i.e., for a universal constant $c > 0$. 
 To justify this equivalence, 
 if, for fixed $\varepsilon > 0$ and for each coloring $\mathcal{C}$, 
 we have that $\liminf_{\ell \to \infty} r_{\mathcal{C}}(\ell) \geq \varepsilon > 0$, 
 then, for every coloring $\mathcal{C}$, we have that $r_{\mathcal{C}}(M) \geq \frac{\varepsilon}{2}$
 for sufficiently large $M$. Conversely, if, for fixed $c > 0$ and for 
 every coloring $\mathcal{C}$, we have that 
 $r_{\mathcal{C}}(M) \geq c > 0$ for sufficiently large $M$, then $\liminf_{\ell \to \infty} r_{\mathcal{C}}(\ell) \geq c > 0$. 

 As indicated above, we prove that no universal constant of the form $c > 0$ indicated above can exist, thus proving 
 \eqref{ESformulation} in the affirmative. 

\begin{example}
 Let each odd and positive integer be colored red. Let each even and positive integer be colored blue. 
 According to this coloring, for a positive integer $M$, let $a(M)$
 denote the number 
 of distinct monochromatic products among $1$, $2$, $\ldots$, $M$. 
 For example, we see that 
 $a(9) = 7$, since the associated monochromatic products are as shown below.

 \

\noindent 1 = {\color{red}1} $\cdot$ {\color{red}1}, \ 
 3 = {\color{red}1} $\cdot$ {\color{red}3}, \ 
 4 = {\color{blue}2} $\cdot$ {\color{blue}2}, \ 
 5 = {\color{red}1} $\cdot$ {\color{red}5}, \ 
 7 = {\color{red}1} $\cdot$ {\color{red}7}, \ 
 8 = {\color{blue}2} $\cdot$ {\color{blue}4}, \ 
 9 = {\color{red}1} $\cdot$ {\color{red}9} = {\color{red}3} $\cdot$ {\color{red}3}.

 \ 

\noindent Observe that 2 = {\color{red}1} $\cdot$ {\color{blue}2} does not qualify, 
 and observe that 
 6 = {\color{red}1} $\cdot$ {\color{blue}6} = {\color{blue}2} $\cdot$ {\color{red}3}
 also does not qualify. It is easily seen (omitting details) that 
 the integer sequence $(a(n) : n \in \mathbb{N})$, 
 as defined above, may be explicitly evaluated
 so that 
\begin{equation}\label{evalan}
 a(n) = \left\lfloor \frac{3 n+1}{4} \right\rfloor. 
\end{equation} 
 This agrees with the sequence indexed in the On-line Encyclopedia of Integer Sequences\footnote{See \url{https://oeis.org/}.}
 as {\tt A037915}, 
 but there is not currently any equivalent definition of $a(n)$, as defined above, 
 given. From the evaluation in \eqref{evalan}, 
 we find that, according to the specified coloring $\mathcal{C}$, 
 and with respect to the notation in \eqref{defineRCN}, 
 the evaluation $\lim_{N \to \infty} r_{\mathcal{C}}(N) = \frac{3}{4}$ holds. 
 \end{example}

 The problem indicated above may be seen in 
 relation to subsequent work related to the given Erd\H os--S\'ark\"ozy references 
 \cite{ErdosSarkozy1990,ErdosSarkozy1993}. 
 Csikv\'{a}ri et al.\ \cite[Theorem 7]{CsikvariGyarmatiSarkozy2012}
 use the logarithmic-density theorem for monochromatic products from 
 Erd\H os and S\'ark\"ozy \cite{ErdosSarkozy1990} 
 to establish positive lower logarithmic density for integers admitting both monochromatic sum and product representations, 
 with the two representations allowed to use different colours. 
  However, this result does not settle the two-class lower-asymptotic-density problem studied here. 
 The Erd\H os--S\'ark\"ozy problem above may also be seen in contrast to a related 
 problem on colorings of sums of integers, as opposed to products of integers, with reference
 to the so-called Erd\H os Problem $\#$484\footnote{See \url{https://www.erdosproblems.com/484}.}. 

 \section{Main construction}
 We begin by covering some basic preliminaries. In this direction, for a subset $S$ of the set $\mathbb{N}$ of positive integers, we write 
 $SS$ in place of $\{ s t : s, t \in S \}$ (where repeated factor are permitted). 
 For subsets $S$ and $T$ contained in $\mathbb{N}$, 
 we write 
\begin{equation}\label{defineP}
 P(S, T) = SS \cup TT. 
\end{equation} 
 
 We require some familiar arithmetic functions, 
 and, in this direction, we 
 recall the prime-counting function 
 $\pi(x)$ giving the number primes $p$ less than or equal to a given real value $x$, 
 and we let 
 \(v_p(n)\) denote the exponent of \(p\) in the prime factorization of the positive integer \(n\). 
 Throughout, sums indexed by \(p\) are understood as being over primes.

 We also require some basic properties of 
 the special function 
 known as the \emph{digamma function}, 
 which may be defined so that 
 $\psi(x) = \frac{\Gamma'(x)}{\Gamma(x)}$ for $x > 0$ and for the familiar $\Gamma$-function, 
 referring to Rainville's classic text for background \cite[\S2]{Rainville1960}. 
 For the Euler--Mascheroni constant $\gamma = \lim_{n \to \infty} \big( 1 + \frac{1}{2} + \cdots + \frac{1}{n} 
 - \log n \big)$, a fundamental property of the digamma function is such that 
\begin{equation}\label{psiexpand} 
 \psi(z) = -\gamma - \frac{1}{z} + \sum_{n=1}^{\infty} \frac{z}{n(z+n)}, 
\end{equation}
 for real $z > 0$. 
 A consequence of a result due to Alzer and Berg 
 \cite[Theorem 4.2(i)]{AlzerBerg2006}
 is such that 
\begin{equation}\label{eqAlzerBerg}
 \psi(b+1) - \psi(a+1) \leq \log\left( \frac{b}{a} \right) 
\end{equation}
 for $0 < a \leq b$, which follows from the mapping 
 $x \mapsto \psi(x+1) - \log x$ being nonincreasing on $(0, \infty)$.

 We require, for the purposes of the proof below, 
 a well known consequence of Markov's inequality. 
 This may be seen as a ``non-probabilistic'' version of Markov's inequality, 
 and this may be formulated as follows. 
 Suppose that $E$ is a finite set. Then, for a mapping $F\colon E \to [0, \infty)$
 and for a positive value $\eta$, 
 the relation 
\begin{equation}\label{displayMarkov}
 \eta \, \big| \{ u \in E : F(u) \geq \eta \} \big| \leq \sum_{u \in E} F(u) 
\end{equation}
 holds.

 Some further notation local to our construction 
 is given as follows. By letting $q \geq 2$ be real, we write
\begin{equation}\label{defineLq} 
 L(q) = \sum_{p \leq q} \frac{1}{p},   
\end{equation}
  noting that the condition $q \geq 2$ guarantees that $L(q) > 0$. 
 Recall that a fundamental result in number theory (attributed to Euler) is such that $L(q) \to \infty$ as $q \to \infty$. 
 Now, for a real value 
 $\theta \in (0, 1)$, define
\begin{equation}\label{defineH}
 H_{\theta}(N, q) 
 = \big| \{ n \in \mathbb{N} : n \leq N, \ \exists \, d\mid n \ \theta q \leq d \leq q \} \big|. 
\end{equation}
 The definition in 
 \eqref{defineH} leads us toward 
 Lemma \ref{lemmaHbound} below.

\begin{lemma}\label{lemmaHbound}
 Let $\theta \in (0, 1)$. Suppose that 
\begin{equation}\label{qgeqmax}
 q \geq \max\{ 2, \theta^{-1} \}, 
\end{equation} 
 and suppose that 
\begin{equation}\label{Ngeqsquare}
 N \geq q^2, 
\end{equation}
 for real values $q$ and $N$. Then
 $$H_{\theta}(N,q) \leq K_{\theta} \frac{N}{L(q)},$$
 for $K_{\theta} = \frac{48}{\theta} + 100 + 48 \log \frac{1}{\theta}$ (noting that $L(q) > 0$ since $q \geq 2$). 
\end{lemma} 

\begin{proof}
 Write $L = L(q)$ for brevity (again noting that $L > 0$ from the hypothesis in \eqref{qgeqmax}). 
 We then define 
\begin{equation}\label{defineFn}
 F(n) = \sum_{p \leq q} \text{{\bf 1}}_{p \mid n}, 
\end{equation}
 letting indicator functions be denoted as per usual. We also define
\begin{equation}\label{defineR}
 R(n) = \sum_{p \leq q} \max\{ v_{p}(n) - 1, 0 \}. 
\end{equation}
 For a real value $ M \geq q$, 
 we find that 
\begin{multline*}
 \sum_{n \leq M} (F(n) - L)^{2}
 = \\ 
 \sum_{p \leq q} \left\lfloor \frac{M}{p} \right\rfloor 
 + 2 \sum_{ p \leq q } \sum_{\substack{p < r \leq q \\ \text{$r$ prime}}}
 \left\lfloor \frac{M}{pr} \right\rfloor 
 - 2 L \sum_{p \leq q} \left\lfloor \frac{M}{p} \right\rfloor 
 + L^2 \lfloor M \rfloor. 
\end{multline*}
 Applying the relation $\lfloor x \rfloor \leq x$ to the positive terms, 
 and applying the relation $\lfloor x \rfloor \geq x -1$ 
 in the negative terms, we find that 
\begin{equation}\label{getridpi}
 \sum_{n \leq M} (F(n) - L)^{2}
 \leq M \left( L - \sum_{p \leq q} \frac{1}{p^2} \right) 
 + 2 L \, \pi(q). 
\end{equation}
 Since 
\begin{equation}\label{comparepiM}
 \pi(q) \leq q \leq M, 
\end{equation} 
 the relation in \eqref{getridpi} gives us that 
\begin{equation}\label{withoutpi} 
 \sum_{n \leq M} (F(n) - L)^{2} \leq 3 M L. 
\end{equation} 
 Now, consider the families among
\begin{equation}\label{familiesamong} 
 \left\{ n \leq M : F(n) < \frac{3L}{4} \right\} 
 \ \text{and} \ \left\{ n \leq M : F(n) \geq \frac{5L}{4} \right\}. 
\end{equation}
 Since $F(n) < \frac{3 L}{4} \Longrightarrow |F(n) - L| > \frac{L}{4}$
 and since $F(n) \geq \frac{5L}{4} \Longrightarrow |F(n) - L| \geq \frac{L}{4}$, 
 either case gives us that 
 $(F(n) - L)^2 \geq \frac{L^2}{16}$. 
 This and \eqref{withoutpi} together give us that 
 each of the families indicated in \eqref{familiesamong}
 has at most $\frac{48 M}{L}$ members.

 Also, by using the evaluation 
 $\sum_{m \geq 2} \frac{1}{m(m-1)} = 1$, we find that 
\begin{equation}\label{boundsumR} 
 \sum_{n \leq M} R(n) = \sum_{p \leq q} \sum_{k \geq 2} \left\lfloor \frac{M}{p^k} \right\rfloor 
 \leq M \sum_{p \leq q} \frac{1}{p(p-1)} \leq M. 
\end{equation}
 Now, since $R(n) \geq 0$ (recalling \eqref{defineR}) and since 
 $L > 0$ (recalling \eqref{defineLq}), 
 by applying the consequence of Markov's inequality
 on display in \eqref{displayMarkov} to 
 the function $R$, with $E = \{ 1, 2, \ldots, \lfloor M \rfloor \}$
 and with $\eta = \frac{L}{4}$, we find that 
\begin{align}
\begin{split}
 \left| \left\{ n \leq M : R(n) > \frac{L}{4} \right\} \right|
 & \leq \left| \left\{ n \leq M : R(n) \geq \frac{L}{4} \right\} \right| \\ 
 & \leq \frac{4}{L} \sum_{n \leq M} R(n), 
\end{split}\label{splitsumR}
\end{align}
 so that a joint application of \eqref{boundsumR} and \eqref{splitsumR} gives us that 
\begin{equation}\label{4MLbound} 
 \left| \left\{ 
 n \leq M : R(n) > \frac{L}{4} \right\} \right| 
 \leq \frac{4 M}{L}. 
\end{equation}

 This leads us toward the following case analysis, where the cases below correspond 
 to the conditions associated with the sets of the forms $D_i$ and $S_d^{(i)}$ 
 defined below.

 \

\noindent\textbf{Case 1.} Define 
\begin{equation}\label{defineD1}
 D_1 =
\left\{
d \in \mathbb{N} :
\theta q \leq d \leq q, \ F(d) < \frac{3L}{4}
\right\}, 
\end{equation}
 and, for each $d \in D_1$, define
\[
S_d^{(1)}
=
\left\{
dk : k \in \mathbb{N}, \ dk \leq N
\right\}.
\]
Taking $M=q$ in the established bound for the first set in
\eqref{familiesamong}, we obtain
\begin{equation}\label{boundcardD1}
|D_1|
\leq
\left|
\left\{
m \in \mathbb{N} :
m \leq q,\quad F(m) < \frac{3L}{4}
\right\}
\right|
\leq \frac{48q}{L}.
\end{equation}
 Again for $d$ in the family defined in \eqref{defineD1}, 
 since $d \geq \theta q$, we find that 
\begin{equation}\label{boundcaardSd1}
 \left|S_d^{(1)}\right|
=
\left\lfloor \frac{N}{d} \right\rfloor
\leq \frac{N}{d}
\leq \frac{N}{\theta q}.
\end{equation}
The union $\bigcup_{d \in D_1} S_d^{(1)}$ is exactly the set of
positive integers $n \leq N$ having at least one divisor in $D_1$.
Consequently, an application of the inclusion--exclusion principle gives us that 
 $\left|\bigcup_{d \in D_1} S_d^{(1)}\right|
 \leq \sum_{d \in D_1} \left|S_d^{(1)}\right|$. 
 This and 
 \eqref{boundcardD1} and 
 \eqref{boundcaardSd1} together give us that 
 \begin{align}
\begin{split}
 \left|\bigcup_{d \in D_1} S_d^{(1)}\right|
 &\leq |D_1|\,\frac{N}{\theta q} \\
 & \leq \frac{48N}{\theta L}.
\end{split}\label{boundcupSd1}
 \end{align}

 \

\noindent {\bf Case 2.} Define 
\begin{equation}\label{defineD2} 
 D_{2} = \left\{ d \in \mathbb{N} : \theta q \leq d \leq q, F(d) \geq \frac{3L}{4} \right\}, 
\end{equation}
 and, for each element $d \in D_{2}$, define 
\begin{equation}\label{defineSd} 
 S_{d}^{(2)} = \left\{ d k : k \in \mathbb{N}, dk \leq N, F(k) < \frac{3L}{4} \right\}. 
\end{equation}
 For a fixed integer $d$ in $D_2$, 
 we find that 
\begin{equation}\label{Noverdgeq} 
 \frac{N}{d} \geq \frac{N}{q}. 
\end{equation} 
 From the hypothesis indicated in \eqref{Ngeqsquare}, 
 we find that \eqref{Noverdgeq} implies that 
\begin{equation}\label{simplifylower} 
 \frac{N}{d} \geq q. 
\end{equation} 
 The relation in \eqref{simplifylower} allows us to let $M = \frac{N}{d}$, recalling that $M \geq q$. 
 So, from the established upper bound for the first family on display in \eqref{familiesamong}, we find that 
\begin{equation}\label{replaceM} 
 \left| \left\{ n \leq \frac{N}{d} : F(n) < \frac{3L}{4} \right\} \right|
 \leq \frac{48 N}{dL}. 
\end{equation} 
 We rewrite \eqref{replaceM} so that 
 $ \big| \big\{ n \leq \frac{N}{d} : F(n) < \frac{3L}{4} \big\} \big|
 \leq \frac{48 N}{dL}$, 
 and this then gives us that 
\begin{equation}\label{rewrittenk} 
 \left| \left\{ k \in \mathbb{N} : k d \leq N, F(k) < \frac{3L}{4} \right\} \right|
 \leq \frac{48 N}{dL}. 
\end{equation}
 For each member $k$ in the family whose cardinality is given 
 as the left-hand side of the inequality in 
 \eqref{rewrittenk}, we obtain a different 
 value of $n = dk \leq N$. 
 So, again for a given \emph{fixed} value $d$, 
 there are (according to \eqref{rewrittenk}) at most $\frac{48 N}{dL}$
 integers $n \leq N$ admitting a factorization of the form $n = dk$
 such that $F(k) < \frac{3L}{4}$. 

 It is immediate that the right-hand side of \eqref{defineSd} has the same cardinality as $ \big\{ k \in \mathbb{N} : k d \leq N, F(k) 
 < \frac{3L}{4} \big\}$. So, from the relation in \eqref{rewrittenk}, we find that $| S_{d}^{(2)} | \leq \frac{48 N}{dL}$. This and an 
 application of the inclusion--exclusion principle give us that 
\begin{align*}
 \left| \bigcup_{d \in D_{2}} S_{d}^{(2)} \right|
 & \leq \sum_{d \in D_{2}} \left| S_{d}^{(2)} \right| \\
 & \leq \frac{48N}{L} \sum_{d \in D_{2}} \frac{1}{d} \\
 & \leq \frac{48N}{L} \sum_{\substack{d \in \mathbb{N} \\ \theta q \leq d \leq q}} \frac{1}{d}. 
\end{align*}
 Since $\psi$ is increasing on $(0, \infty)$ (and this may be seen via the series expansion shown in \eqref{psiexpand}), 
 and since the recurrence $\psi(x+1) - \psi(x) = \frac{1}{x}$ holds, 
 a telescoping argument and the Alzer--Berg 
 inequality in \eqref{eqAlzerBerg} (applied with $a = \theta q$ and $b = q$)
 give us that 
\begin{align*} 
 \sum_{\substack{d \in \mathbb{N} \\ \theta q \leq d \leq q } } \frac{1}{d} 
 & = \sum_{\substack{d \in \mathbb{N} \\ \theta q \leq d \leq q}} 
 \big( \psi(d+1) - \psi(d) \big) \\ 
 & \leq \psi(q+1) - \psi(\theta q) \\ 
 & = \frac{1}{\theta q} + \psi(q+1) - \psi(\theta q +1) \\
 & \leq \frac{1}{\theta q} + \log\left( \frac{1}{\theta} \right) \\ 
 & \leq 1 + \log\left( \frac{1}{\theta} \right), 
\end{align*}
 noting that this last inequality holds since $q \geq \theta^{-1}$, 
 i.e., from the hypothesis displayed in \eqref{qgeqmax}. 
 A combined application of the given bounds for $ \left| \bigcup_{d \in D_{2}} S_{d}^{(2)} \right|$ 
 and $ \sum_{\substack{d \in \mathbb{N} \\ \theta q \leq d \leq q } } \frac{1}{d}$ gives us that 
\begin{equation}\label{boundcupSd2}
 \left| \bigcup_{d \in D_{2}} S_{d}^{(2)} \right| 
 \leq \frac{48 N}{L} \left( 1 + \log\left( \frac{1}{\theta} \right) \right). 
\end{equation}

 \ 

\noindent {\bf Case 3.} Let $D_2$ be as in \eqref{defineD2}, and, for each member $d$ in $D_2$, define 
\begin{equation}\label{defineSd3}
 S_{d}^{(3)} = \left\{ d k : k \in \mathbb{N}, \ dk \leq N, \ F(k) \geq \frac{3L}{4} \right\}. 
\end{equation}
 Let $n \in S_d^{(3)}$, and write $k = \frac{n}{d}$. 
 From the definition in \eqref{defineFn}, we find that 
\begin{equation}\label{combinationF}
 F(d) + F(k) - F(n) = \sum_{\substack{p \leq q \\ p\mid d, p \mid k }} 1. 
\end{equation}
 Now, observe that a given prime $p \leq q$ dividing both $d$ and $k$ is 
 such that $v_{p}(n) \geq 2$ and thus contributes 
 at least $1$ as a term on the right-hand side of \eqref{defineR}. 
 From \eqref{combinationF}, we thus have that 
\begin{equation}\label{compareFR} 
 F(d) + F(k) - F(n) \leq R(n). 
\end{equation}
 From \eqref{compareFR} and the assumptions that $F(d) \geq \frac{3L}{4}$ 
 and that $F(k) \geq \frac{3L}{4}$, we find that 
\begin{align}
\begin{split}
 F(n) & \geq F(d) + F(k) - R(n) \\ 
 & \geq \frac{3L}{2} - R(n). 
\end{split}\label{splitFLR}
\end{align}
 From \eqref{splitFLR}, we obtain the implication $R(n) \leq \frac{L}{4} \Longrightarrow F(n) \geq \frac{5L}{4}$. So, each member $n$ in 
 $S_{d}^{(3)}$ is such that $R(n) > \frac{L}{4}$ or $F(n) \geq \frac{5L}{4}$. 
 Since this holds for an arbitrary member $d$ of $D_{2}$, 
 and since the definition in \eqref{defineSd3}
 gives us that each member in $S_{d}^{(3)}$
 is $\leq N$, we find that the inclusion
\begin{multline*}
 \bigcup_{d \in D_{2}} S_{d}^{(3)} \subseteq \\ 
 \left\{ n \in \mathbb{N} : n \leq N, \ R(n) > \frac{L}{4} \right\} 
 \cup \left\{ n \in \mathbb{N} : n \leq N, \ F(n) \geq \frac{5L}{4} \right\} 
\end{multline*}
 holds, and we thus obtain that 
\begin{multline*}
 \left| \bigcup_{d \in D_{2}} S_{d}^{(3)} \right| \leq \\ 
 \left| \left\{ n \in \mathbb{N} : n \leq N, \ R(n) > \frac{L}{4} \right\} \right| 
 + \left| \left\{ n \in \mathbb{N} : n \leq N, \ F(n) \geq \frac{5L}{4} \right\} \right|. 
\end{multline*}

 From \eqref{qgeqmax} and \eqref{Ngeqsquare} together, we find that 
 $N \geq q$. This allows us (recalling that $M \geq q$, as in \eqref{comparepiM}) to set $M = N$, 
 i.e., in the established bounds for 
 the latter family in \eqref{familiesamong} 
 and in the bound on display in \eqref{4MLbound}. 
 Using these established bounds, 
 the above bound for $ \left| \bigcup_{d \in D_{2}} S_{d}^{(3)} \right| $ gives us that 
\begin{align}
\begin{split}
 \left| \bigcup_{d \in D_{2}} S_{d}^{(3)} \right|
 & \leq \frac{4N}{L} + \frac{48N}{L} \\
 & = \frac{52N}{L}. 
\end{split}\label{boundcupSd3}
\end{align}

 Now, recall the definition of 
 $H_{\theta}(N, q)$ on display in \eqref{defineH}, 
 giving a family of positive integers $n \leq N$ admitting a factorization
 $n = k d$ for positive integers $k$ and $d$, with $ \theta q \leq d \leq q$. 
 Each such factorization satisfies one of the conditions used to define Cases 1--3, so that 
\begin{equation}\label{splitHcases}
 H_{\theta}(N, q) \leq \left| \bigcup_{d \in D_{1}} S_{d}^{(1)} \right| 
 + \left| \bigcup_{d \in D_{2}} S_{d}^{(2)} \right| + \left| \bigcup_{d \in D_{2}} S_{d}^{(3)} \right|. 
\end{equation}
 So, a combined application of 
 the inequalities among \eqref{boundcupSd1}, 
 \eqref{boundcupSd2}, \eqref{boundcupSd3}, and \eqref{splitHcases}
 gives us that 
\begin{align*}
 H_{\theta}(N, q) 
 & \leq \frac{48 N}{\theta L} + \frac{48 N}{L} \left( 1 + \log \frac{1}{\theta} \right)
 + \frac{52 N}{L} \\ 
 & = K_{\theta} \frac{N}{L(q)}, 
\end{align*}
 as desired. 
\end{proof}

\begin{lemma}\label{fourpart}
 Let $\theta \in (0, 1)$ and let $\delta > 0$. Then there exist a strictly increasing sequence $(N_j : j \in \mathbb{N})$ of positive integers 
 and a set 
\begin{equation}\label{Sirrational}
 S \subseteq \big( \max\{ 2, \theta^{-1} \}, \infty \big) \setminus \mathbb{Q} 
\end{equation}
 such that: 

\begin{enumerate}[label=(\alph*)]

\item Each bounded interval contains only finitely many elements of $S$; 

\item For each positive integer $j$, the value $\sqrt{N_j}$ is in $S$; 

\item The elements in $S$ are such that 
 $$ \sum_{s \in S} \frac{1}{L(s)} \leq \delta\text{; and} $$

\item The function 
\begin{equation}\label{definechi}
 \chi(x) = (-1)^{| S \cap (1, x) |} 
\end{equation}
 defined for $x \in \mathbb{Q}$ and $x \geq 1$ is such that 
\begin{equation}\label{chireflect} 
 \chi(x) = -\chi\left( \frac{N_j}{x} \right) 
\end{equation}
 for $j \geq 1$ and 
\begin{equation}\label{conditionreflect}
 x \in [1, N_j] \cap \mathbb{Q}.
\end{equation}

\end{enumerate}

\end{lemma}

\begin{proof}
 Our strategy is to construct the terms of the desired sequence $(N_j : j \in \mathbb{N})$ recursively, starting with $N_{0} = 1$. Similarly, 
 we recursively construct a sequence $(S_j : j \in \mathbb{N})$ of finite sets to define $S$, with 
\begin{equation}\label{Sbigcup} 
 S = \bigcup_{j \geq 1} S_{j}. 
\end{equation} 
 Also, define $S_{0} = \varnothing$. Now, for an index $j \geq 1$, 
 with the assumption that $N_{j-1}$ and $S_{j-1}$ have already been determined, 
 we proceed to choose a positive integer $M_{j}$ satisfying both 
\begin{equation}\label{firstMcondition} 
 M_{j} > \max\{ N_{j-1}, 2, \theta^{-1} \} 
\end{equation}
 and 
\begin{equation}\label{definingLMj} 
 L(M_{j}) \geq \frac{4^{j}}{2 \delta}, 
\end{equation}
 noting that such a selection is permitted by the 
 classical result giving that the right-hand side of \eqref{defineLq} approaches 
 $\infty$ as $q \to \infty$.

 Having chosen $M_j$ (subject to the given conditions), we proceed to define
\begin{equation}\label{defineNj}
 N_{j} = M_{j}^{2} + 1
\end{equation}
 and 
\begin{equation}\label{Sjrecursion} 
 S_{j} = S_{j-1} \cup \left\{ \sqrt{N_{j}} \right\} \cup \left\{ \frac{N_{j}}{s} : s \in S_{j-1} \right\}. 
\end{equation}
 From \eqref{firstMcondition} and \eqref{defineNj} together, we obtain that 
\begin{equation}\label{Nstrict} 
 N_{j} = M_{j}^2 + 1 > M_j > N_{j-1}
\end{equation}
 for $j \geq 1$, thus giving us that 
 $(N_{j} : j \in \mathbb{N})$ 
 is strictly increasing and consists of positive integers. 

 Now, since $M_j$ is a positive integer, and since $M_j^2 < N_j = M_j^2 + 1 < (M_j + 1)^2$, we find that $N_j$ is not a perfect square, so 
 that $\sqrt{N_j}$ is irrational for each index $j \geq 1$. 
 Since $S_0 = \varnothing$, the recursion in \eqref{Sjrecursion}
 gives us that $S_{1} = \left\{ \sqrt{N_{1}} \right\}$. 
 Since $\sqrt{N_1}$ is positive and irrational (as established), 
 we find that $S_1$ is a finite set of positive irrational numbers. 
 Also observe that $\frac{N_j}{s}$ is positive and irrational whenever $s$ is positive and irrational. 
 So, the recursive definition for $S_j$ and the irrationality of 
 $\sqrt{N_j}$ together show that $S_j$ is necessarily a finite set consisting of positive and irrational numbers, 
 for each index $j \geq 1$. 
 We proceed to prove, via induction, that 
\begin{equation}\label{inductinclude} 
 S_{j} \subseteq \big( \max\{ 2, \theta^{-1} \}, N_{j} \big) 
\end{equation}
 for each index $j \geq 1$.

 For the base case, with $j = 1$, 
 we find that $S_1 = \left\{ \sqrt{N_1} \right\}$ 
 and that $\max\big\{ 2, \theta^{-1} \big\} < M_1$, 
 directly from the definition in \eqref{firstMcondition}. 
 Similarly, from \eqref{defineNj}, 
 we find that $M_1 < \sqrt{N_1} < N_1$. 
 So, since $\max\big\{ 2, \theta^{-1} \big\} < \sqrt{N_1} < N_1$, 
 the inclusion in \eqref{inductinclude} holds for the current base case. 
 Now, for an index $j \geq 2$, we (inductively) assume that 
\begin{equation}\label{inductS}
 S_{j-1} \subseteq \big( \max\{ 2, \theta^{-1} \}, N_{j-1} \big). 
\end{equation}
 Now, for each element $s \in S_{j-1}$, 
 the inductive hypothesis in \eqref{inductS}
 gives us that 
 $ \max\{ 2, \theta^{-1} \} < s < N_{j-1}$. 
 Also, from \eqref{firstMcondition}, we find that 
\begin{equation}\label{NstrictM}
 N_{j-1} < M_j. 
\end{equation}
 Also, from \eqref{defineNj}, we find that 
\begin{equation}\label{MsqrtNN} 
 M_j < \sqrt{N_j} < N_j, 
\end{equation}
 i.e., so that 
\begin{equation}\label{firstunion} 
 \max\{ 2, \theta^{-1} \} < s < N_{j}. 
\end{equation}
 Again letting $s \in S_{j-1}$, from the induction hypothesis shown in \eqref{inductS}, we find that $s < N_{j-1}$. 
 Moreover, from \eqref{NstrictM} 
 and \eqref{MsqrtNN}, 
 we find that 
\begin{equation}\label{sNMsqrtN} 
 s < N_{j-1} < M_j < \sqrt{N_j}. 
\end{equation}
 From \eqref{sNMsqrtN}, we see that $1 < \max\big\{ 2, \theta^{-1} \big\} < s < \sqrt{N_j}$, 
 which, in turn, implies 
\begin{equation}\label{secondunion}
 \max\big\{ 2, \theta^{-1} \big\} 
 < \sqrt{N_j} < \frac{N_j}{s} < N_j. 
\end{equation}
 From \eqref{firstunion} and \eqref{secondunion}, 
 we obtain the inclusion
\begin{align}
\begin{split}
 S_{j} & = S_{j-1} \cup \left\{ \sqrt{N_{j}} \right\} \cup \left\{ \frac{N_{j}}{s} : s \in S_{j-1} \right\} \\ 
 & \subseteq \big( \max\{ 2, \theta^{-1} \}, N_{j} \big), 
\end{split}\label{splitinclusion}
\end{align}
 thus completing our induction. 

 Now, recall that we previously established that 
 $S_j$ consists of irrational numbers. 
 From this and \eqref{splitinclusion} together, we find that 
 $$ S = \bigcup_{j \geq 1} S_{j} \subseteq \big( \max\{ 2, \theta^{-1} \}, \infty \big) \setminus \mathbb{Q}. $$

 To prove the first assertion listed above, we begin by proving that 
\begin{equation}\label{diffinclusion} 
 S_{j} \setminus S_{j-1} \subseteq (N_{j - 1}, \infty) 
\end{equation}
 for each $j \geq 1$. Now, recall the relations 
 among $N_{j-1}$, $M_j$, and $\sqrt{N_j}$ shown in \eqref{sNMsqrtN}. By letting $s \in S_{j-1}$, we find, 
 again with reference to \eqref{sNMsqrtN} (recalling that 
 \eqref{sNMsqrtN} holds for $s \in S_{j-1}$), 
 that $0 < s < \sqrt{N_j}$.
 So, the relations in \eqref{sNMsqrtN} give us that 
\begin{equation}\label{furtherrelations}
 \frac{N_{j}}{s} > \frac{N_j}{\sqrt{N_j}} = \sqrt{N_j} > N_{j-1}. 
\end{equation}
 Now, from the recursion on display in \eqref{Sjrecursion}, if $t$ is any element in $S_{j} \setminus S_{j-1}$, 
 then either $t = \sqrt{N_{j}}$ or $t = \frac{N_{j}}{s}$ for some $s \in S_{j-1}$. 

In the former case, the definitions of $N_{j}$ and $M_{j}$ give us that 
 $$ t = \sqrt{N_{j}} = \sqrt{M_{j}^{2} + 1} > M_{j} > N_{j-1}. $$
 In the latter case, we have $t = \frac{N_{j}}{s}$ for some $s \in S_{j-1}$, 
 and the previously established inequality $s < \sqrt{N_{j}}$ gives us that 
 $$ t = \frac{N_{j}}{s} > \sqrt{N_j} > N_{j-1}. $$
 So, for every $j \geq 1$, the inclusion
 $$ S_{j} \setminus S_{j-1} \subseteq (N_{j-1}, \infty) $$
 holds.

 We proceed to let 
 $I \subseteq \mathbb{R}$ be an arbitrary bounded interval, 
 and we then let $T $ be a real value such that 
 $I \subseteq (-\infty, T]$. 
 As established previously (recalling \eqref{Nstrict}), 
 the sequence $(N_{j})_{j \geq 1}$ is strictly increasing, and we also recall that it consists of positive integers. 
 So, since $N_{j} \to \infty$, 
 we are permitted to choose an index $J \geq 1$
 such that $N_J > T$. 
 For each $j > J$, every member $t$ in $S_{j} \setminus S_{j-1}$ is such that 
 $$ t > N_{j-1} \geq N_{J} > T, $$
 and we thus find that 
\begin{equation}\label{minuscapnull}
 \big( S_{j} \setminus S_{j-1} \big) \cap I = \varnothing, 
\end{equation} 
 since $I \subseteq (-\infty, T]$. The recursion in \eqref{Sjrecursion} gives us that $S_{j-1} \subseteq S_{j}$, and thus, for $j > J$, we 
 proceed to write
\begin{equation}\label{simplifynull}
 S_{j} \cap I = (S_{j-1} \cap I) \cup \big( (S_{j} \setminus S_{j-1}) \cap I \big), 
\end{equation}
 so that a combined application of \eqref{minuscapnull} and \eqref{simplifynull} yields
\begin{equation}\label{userepeatedly} 
 S_{j} \cap I = S_{j-1} \cap I. 
\end{equation}
 A repeated application of the equality in \eqref{userepeatedly} then gives us that 
\begin{equation}\label{SupperJcapI}
 S_{j} \cap I = S_{J} \cap I 
\end{equation}
 for each $j \geq J$. 
 Since 
\begin{equation}\label{chainS}
 S_{j} \subseteq S_{J} 
 \ \text{for} \ j \in \{ 1, 2, \ldots, J-1 \}, 
\end{equation}
 we find, by appealing to 
 the equality in \eqref{SupperJcapI} that 
\begin{align*}
 S \cap I 
 & = \bigcup_{j \geq 1} \left( S_{j} \cap I \right) \\ 
 & = S_{J} \cap I. 
\end{align*}
 Also observe that $S_{J} \cap I$ is finite, as a consequence of $S_{J}$ being finite. This proves the first listed assertion above. 

 The second listed assertion follows immediately from the recursion in \eqref{Sjrecursion} in conjunction with the definition 
 in \eqref{Sbigcup}. 

 We now proceed to establishing the third assertion listed above. In this direction, we begin by observing that $S_{j-1}$, $\left\{ 
 \sqrt{N_{j}} \right\}$, and $\big\{ \frac{N_{j}}{s} : s \in S_{j-1} \big\}$ are pairwise disjoint for an arbitrary element $j$ in $\mathbb{N}$, 
 since: For each $s \in S_{j-1}$, we have that 
 $s < \sqrt{N_j} < \frac{N_j}{s}$, recalling 
 \eqref{sNMsqrtN} and \eqref{furtherrelations}. 
 For a fixed and positive integer $j$, we find that the mapping on $S_{j-1}$
 such that $s \mapsto \frac{N_{j}}{s}$ is injective, 
 and we thus find that the recursion shown in 
 \eqref{Sjrecursion} yields the recurrence such that 
\begin{equation}\label{cardrec} 
 |S_{j}| = 2 |S_{j-1}| + 1. 
\end{equation}
 Since $S_{0} = \varnothing$, this and the recurrence in \eqref{cardrec} together give us that 
\begin{equation}\label{cardeval} 
 |S_{j}| = 2^{j} - 1
\end{equation}
 for each nonnegative integer $j$. 
 Since $S_{j-1} \subseteq S_{j}$ (again recalling the recursion in \eqref{Sjrecursion}), 
 this and the evaluation in \eqref{cardeval} together give us that 
\begin{equation}\label{carddiff}
 \big| S_{j} \setminus S_{j-1} \big| = |S_{j}| - |S_{j-1}| = 2^{j-1}
\end{equation}
 for each positive integer $j$. From the recursion in \eqref{Sjrecursion}, 
 a given element 
 $t$ in $S_{j} \setminus S_{j-1}$ (letting $j \geq 1$ be fixed) is either equal to $\sqrt{N_{j}}$ or is equal to $\frac{N_{j}}{s}$
 for some $s \in S_{j-1}$. 
 In the latter case, the relations in \eqref{furtherrelations} (which hold, as above, for $s \in S_{j-1}$)
 give us that $t = \frac{N_{j}}{s} > \sqrt{N_{j}}$, so, in either case, 
 we see that $t \geq \sqrt{N_{j}}$. 
 Now, since $t \geq \sqrt{N_j} > M_j$ (recalling \eqref{MsqrtNN}), 
 and since the mapping on $\mathbb{R}_{\geq 2}$ such that $q \mapsto L(q)$
 is nondecreasing (recalling \eqref{defineLq}), the defining condition
 in \eqref{definingLMj} gives us that 
\begin{equation}\label{bound1overL} 
 \frac{1}{L(t)} \leq \frac{1}{L(M_j)} \leq \frac{2\delta}{4^{j}}. 
\end{equation}
 From \eqref{carddiff} and \eqref{bound1overL} together, we find that 
\begin{align}
\begin{split} 
 \sum_{t \in S_{j} \setminus S_{j-1}} \frac{1}{L(t)}
 & \leq 2^{j-1} \frac{2\delta}{4^{j}} \\
 & = \frac{\delta}{2^{j}}.
\end{split}\label{deltaoverpower}
\end{align}
 Since $S = \bigcup_{j \geq 1} S_{j}$, $S_0 = \varnothing$, and $S_{j-1} \subseteq S_{j}$ for each $j \geq 1$, 
 sets of the form $S_{j} \setminus S_{j-1}$ are pairwise disjoint with union $S$, an application
 of \eqref{deltaoverpower} gives us that 
\begin{align*}
 \sum_{s \in S} \frac{1}{L(s)}
 & = \sum_{j=1}^{\infty} \sum_{t \in S_{j} \setminus S_{j-1}} \frac{1}{L(t)} \\ 
 & \leq \sum_{j=1}^{\infty} \frac{\delta}{2^{j}} \\
 & = \delta,
\end{align*}
 thus proving the third listed assertion as desired.

 We proceed to prove the final listed assertion. In this direction, we 
 begin by fixing a positive integer $j$, and observing that 
 \eqref{inductinclude} gives us that 
\begin{equation}\label{Sjinclude}
 S_{j} \subseteq (1, N_{j}). 
\end{equation} 
 Consequently, the relation 
\begin{equation}\label{Scap1NjSj} 
 S \cap (1, N_j) = S_j 
\end{equation}
 holds. From the chain of sets indicated in \eqref{chainS}
 together with the inclusion in \eqref{diffinclusion}
 together with the $N$-sequence being nondecreasing, 
 we find that 
\begin{align*}
 S \setminus S_j 
 & = \bigcup_{k = j + 1}^{\infty} (S_k \setminus S_{k-1}) \\
 & \subseteq (N_j, \infty).
\end{align*}
 Now, observe that the mapping on $S_j$ such that $t \mapsto \frac{N_{j}}{t}$ is a bijection, since, in view of the recursion in 
 \eqref{Sjrecursion}, the given mapping interchanges $S_{j-1}$ and $\big\{ \frac{N_{j}}{s} : s \in S_{j-1} \big\}$
 and fixes $\sqrt{N_j}$. 
 Now, let 
\begin{equation}\label{xin1Nj} 
 x \in [1, N_j] \cap \mathbb{Q}
\end{equation}
 be fixed. 
 Since $x$ and $\frac{N_{j}}{x}$ are rational and since each element in $S_j$ is irrational, 
 we find that neither $x$ nor $\frac{N_{j}}{x}$ is in $S_{j}$. 
 Now, since the given bijection maps
 $S_j$ onto itself, we see that this same bijection maps
 $S_j \cap (1, x)$ onto $S_j \cap \big( \frac{N_{j}}{x}, N_{j} \big)$, 
 and we thus obtain the equality 
\begin{equation}\label{Sjcap1x} 
 |S_{j} \cap (1, x)| = \left| S_{j} \cap \left( \frac{N_{j}}{x}, N_{j} \right) \right|. 
\end{equation}
 Recalling the inclusion in \eqref{Sjinclude}, 
 since $\frac{N_{j}}{x} \not\in S_j$, 
 we find that 
 $S_j \cap \big( 1, \frac{N_{j}}{x} \big)$ and $S_{j} \cap \big( \frac{N_{j}}{x}, N_{j} \big)$
 form a partition of $S_{j}$. 
 This partition gives us that 
\begin{equation}\label{Sj2set} 
 |S_{j}| = \left| S_{j} \cap \left( 1, \frac{N_{j}}{x} \right) \right| 
 + \left| S_{j} \cap \left( \frac{N_{j}}{x}, N_{j} \right) \right|. 
\end{equation}
 A combined application of \eqref{cardeval}, \eqref{Sjcap1x}, 
 and \eqref{Sj2set} then gives us that 
\begin{align}
\begin{split}
 \left| S_{j} \cap (1, x) \right| + \left| S_{j} \cap \left( 1, \frac{N_{j}}{x} \right) \right|
 & = |S_{j}| \\ 
 & = 2^{j} - 1.
\end{split}\label{equalscardSj}
\end{align}
 Recalling \eqref{xin1Nj}, we have that $x \leq N_j$, with $(1, x) \subseteq (1, N_j)$, so that 
\begin{align}
\begin{split} 
 S \cap (1, x) & = \big( S \cap (1, N_j) \big) \cap (1, x) \\ 
 & = S_{j} \cap (1, x), 
\end{split}\label{associativecap}
\end{align}
 where the last equality in \eqref{associativecap} follows from \eqref{Scap1NjSj}. From \eqref{xin1Nj}, we find that 
 $ 1 \leq \frac{N_{j}}{x} \leq N_{j}, $ 
 i.e., so that $\big( 1, \frac{N_{j}}{x} \big) \subseteq (1, N_{j})$. 
 By analogy with \eqref{associativecap}, we have that 
\begin{align}
\begin{split}
 S \cap \left( 1, \frac{N_{j}}{x} \right)
 & = \big( S \cap (1, N_j) \big) \cap \left( 1, \frac{N_{j}}{x} \right) \\ 
 & = S_{j} \cap \left( 1, \frac{N_{j}}{x} \right), 
\end{split}\label{furtherassociative}
\end{align} 
 where we have again exploited \eqref{Scap1NjSj}. 
 So, a combined application of 
 \eqref{equalscardSj}, \eqref{associativecap}, 
 and \eqref{furtherassociative} yields
\begin{align}
\begin{split} 
 \left| S \cap (1, x) \right| + \left| S \cap \left( 1, \frac{N_{j}}{x} \right) \right| 
 & = \left| S_{j} \cap (1, x) \right| + \left| S_{j} \cap \left( 1, \frac{N_{j}}{x} \right) \right| \\ 
 & = 2^{j} - 1.
\end{split}\label{2tojminus1}
\end{align}
 So, in view of the parity of $2^{j} - 1$, the definition of $\chi$ shown in 
 \eqref{definechi} and the equalities in 
 \eqref{2tojminus1} together give us that 
\begin{align*}
 \chi(x) \chi\left( \frac{N_{j}}{x} \right)
 & = (-1)^{ \left| S \cap (1, x) \right| + \left| S \cap \left( 1, \frac{N_{j}}{x} \right) \right| } \\ 
 & = (-1)^{2^{j} - 1} \\ 
 & = -1.
\end{align*}
 Since the possible values for $\chi$ are in $\{ -1, 1 \}$, the preceding equalities give us that 
 $\chi(x) = -\chi\left( \frac{N_{j}}{x} \right)$, 
 and this proves the final listed assertion. 
\end{proof}

 We are now in a position to prove our main result and to provide an affirmative solution to the Erd\H{o}s--S\'{a}rk\"{o}zy problem given 
 in Section \ref{secintro}. 

\begin{theorem}
 For each $\varepsilon > 0$, there exists a partition $\{ \mathcal{A}_{1}, \mathcal{A}_{2} \}$
 of $\mathbb{N}$ such that 
 $$ \liminf_{N \to \infty} 
 \frac{ \left| P(\mathcal{A}_{1}, \mathcal{A}_{2}) \cap \{ 1, 2, \ldots, N \} \right| }{N} < \varepsilon. $$
\end{theorem}

\begin{proof} 
 It suffices to prove the desired result for $\varepsilon \in (0, 1)$ (since, once such a result is established for $\varepsilon \in (0, 1)$, it would immediately hold for larger values). So, let $\varepsilon \in (0, 1)$ be fixed. Also, 
 set $\theta = \frac{\varepsilon}{2}$ and $\delta = \frac{\varepsilon}{4 K_{\theta}}$. 
 Since $\theta \in (0, 1)$ and $\delta > 0$, Lemma \ref{fourpart}
 gives us a strictly increasing sequence $(N_{j} : j \in \mathbb{N})$ of positive integers and a set $S$ 
 satisfying the desired conditions in 
 Lemma \ref{fourpart}, along
 with the associated function 
 $\chi(x) = (-1)^{|S \cap (1, x)|}$
 for $x \in \mathbb{Q}_{ \geq 1}$. 
 We then define two color classes so that 
 $\mathcal{A}_{1} = \{ n \in \mathbb{N} : \chi(n) = 1 \}$
 and $\mathcal{A}_{2} = \{ n \in \mathbb{N} : \chi(n) =-1 \}$. 
 Since $\chi(n) \in \{ -1, 1 \}$ for each positive integer $n$, 
 we find that $\mathcal{A}_{1}$ and $\mathcal{A}_{2}$ form a partition of $\mathbb{N}$, 
 defining a $2$-coloring that we denote with $\mathcal{C}$. 

 Now, let $j$ be a fixed and positive integer, and write $N$ in place of $N_j$. 
 For any $n \in P(\mathcal{A}_{1}, \mathcal{A}_{2})$ such that 
\begin{equation}\label{thetaNnleqN} 
 \theta N < n \leq N, 
\end{equation} 
 there exists (recalling \eqref{defineP}) 
 at least one factorization $n = a b$ such that $a, b \in \mathbb{N}$ and $a \leq b$ and
\begin{equation}\label{chiachib} 
 \chi(a) = \chi(b). 
\end{equation} 
 We fix one such factorization $n = ab$. 
 This leads us to two cases covered below. 

 \ 

\noindent {\bf Case 1.}
 Suppose that $b \leq \sqrt{N}$. Then, by taking $s = \sqrt{N} \in S$, as permitted by the second 
 listed assertion in Lemma \ref{fourpart}, we find that 
 $ \theta s = \theta \sqrt{N} \leq \frac{\theta N}{b}$. 
 From \eqref{thetaNnleqN}, we 
 find that $\frac{\theta N}{b} < \frac{n}{b} = a \leq b \leq \sqrt{N} = s$.

 \

\noindent {\bf Case 2.} Suppose that $b > \sqrt{N}$. Define 
\begin{equation}\label{uNbinQ} 
 u = \frac{N}{b} \in \mathbb{Q}.
\end{equation} 
 Since $\theta N < n = ab \leq N$ and since $b > \sqrt{N}$, we find that 
\begin{equation}\label{case2string} 
 1 \leq a = \frac{n}{b} \leq \frac{N}{b} = u < \sqrt{N},
\end{equation}
 and 
\begin{equation}\label{miscstring}
 a = \frac{n}{b} > \frac{\theta N}{b} = \theta u. 
\end{equation}
 From \eqref{case2string}, we find that $u \in [1, N]$, 
 and this and \eqref{uNbinQ} together give us that 
 $u \in [1, N] \cap \mathbb{Q}$, 
 recalling the condition in \eqref{conditionreflect} associated with the reflection identity
 in \eqref{chireflect} for the $\chi$-function. 
 This same reflection identity, applied with 
 $u \in [1, N] \cap \mathbb{Q}$, gives us that 
\begin{align*} 
 \chi(u) 
 & = -\chi\left( \frac{N}{u} \right) \\ 
 & = -\chi\left( b \right), \\
 & = -\chi\left( a \right). 
\end{align*} 
 where we have exploited both \eqref{chiachib} and \eqref{uNbinQ}. This gives us that 
\begin{equation}\label{differentcolors} 
 \chi(a) \neq \chi(u). 
\end{equation} 
 Now, since $n = ab \leq N$ and $b > 0$, we have that 
\begin{equation}\label{alequ} 
 a = \frac{n}{b} \leq \frac{N}{b} = u.
\end{equation} 
 Also, observe that \eqref{differentcolors} and \eqref{alequ} together give
 us that $a < u$. 
 Each member of $S$ is irrational, recalling \eqref{Sirrational}, 
 but $a$ is an integer and $u$ is rational. This leads us to the implications such that 
\begin{align}
\begin{split}
 S \cap (a, u) = \varnothing & \Longrightarrow
 S \cap (1, a) = S \cap (1, u) \\ 
 & \Longrightarrow \chi(a) = \chi(u), 
\end{split}\label{implyseq}
\end{align}
 thus contradicting \eqref{differentcolors}, noting that we have exploited (to obtain the final implication in \eqref{implyseq})
 the definition of $\chi$ on display in \eqref{definechi}. 
 So, in view of the implications shown in \eqref{implyseq}, 
 it cannot be the case that $S \cap (a, u)$ is empty. 
 So, we have shown that there exists an element $s \in S$ (noting that $S$ is nonempty by construction, with, for example, 
 the value $\sqrt{N_{\ell}}$ being in $S$ for all $\ell \geq 1$) such that 
 $a < s < u$. 
 For the given choice of $s$, the inequality
 $u < \sqrt{N}$ (recalling \eqref{case2string})
 and the inequality $a > \theta u$ (recalling \eqref{miscstring}) 
 give us that 
 $s < u < \sqrt{N}$ and that $\theta s < \theta u < a < s$. 

 \ 

 So, in either of the above cases, 
 we have that $s \leq \sqrt{N}$ and that $\theta s \leq a \leq s$.

 Now, let $m \in P(\mathcal{A}_{1}, \mathcal{A}_{2}) \cap [1, N]$. 
 So, in other words, the integer $m$ is a monochromatic product such that $1 \leq m \leq N$. 
 We proceed to consider two separate cases. 

 \ 

\noindent {\bf Case 1.} Suppose that $m \leq \theta N$. 
 Then it is immediate that 
 $ m \in \{ k \in \mathbb{N} : k \leq \theta N \}$. 

 \ 

\noindent {\bf Case 2.} Suppose that $\theta N < m \leq N$. From the preceding case analysis, 
 we deduce that there exists a positive integer divisor 
 $a \mid m$ and an element $s \in S$
 such that $s \leq \sqrt{N}$ and $\theta s \leq a \leq s$. 

 \ 

 We thus obtain 
 the inclusion 
\begin{multline*}
 P(\mathcal{A}_{1}, \mathcal{A}_{2}) \cap [1, N]
 \subseteq \{ m \in \mathbb{N} : m \leq \theta N \} \cup \\ 
 \bigcup_{\substack{t \in S \\ t \leq \sqrt{N}}}
 \big\{ m \in \mathbb{N} : m \leq N, \ \exists d \in \mathbb{N} \ 
 \big( d \mid m \ \text{and} \ \theta t \leq d \leq t \big) \big\}. 
\end{multline*}

 Since $S \subset (1, \infty)$ (recalling \eqref{Sirrational}), we find that indexing sets associated with the latter union above consist of 
 elements $t$ in $S \cap [1, \sqrt{N}]$. As established in Lemma \ref{fourpart}, 
 every bounded interval contains only finitely many elements of $S$, 
 and hence the finitude of $S \cap [1, \sqrt{N}]$. 
 From this finitude, the above inclusion, 
 and the definition in \eqref{defineH}, 
 we find that 
\begin{equation}\label{sumH} 
 \big| P(\mathcal{A}_{1}, \mathcal{A}_{2}) \cap [1, N] \big| 
 \leq \lfloor \theta N \rfloor
 + \sum_{\substack{ t \in S \\ t \leq \sqrt{N} }} H_{\theta}(N, t). 
\end{equation}

 Now, for $t \in S$ such that $t \leq \sqrt{N}$, since a given element in $S$ is strictly greater than $\max\{ 2, \theta^{-1} \}$ according to 
 \eqref{Sirrational}, we find that $t > \max\{ 2, \theta^{-1} \}$. Since $t \leq \sqrt{N}$, 
 this yields $N \geq t^2$. 
 So, the required coinditions of Lemma \ref{lemmaHbound} hold for $q = t$, 
 writing 
\begin{equation}\label{repeatlemma} 
 H_{\theta}(N, t) \leq K_{\theta} \frac{N}{L(t)}. 
\end{equation}
 From \eqref{sumH} and \eqref{repeatlemma}, we find that 
\begin{align}
\begin{split}
 \frac{ \left| P(\mathcal{A}_{1}, \mathcal{A}_{2}) \cap [1, N] \right| }{N}
 & \leq \frac{ \lfloor \theta N \rfloor }{N} + K_{\theta} \sum_{\substack{t \in S \\ t \leq \sqrt{N}}} \frac{1}{L(t)} \\ 
 & \leq \theta + K_{\theta} \sum_{\substack{t \in S \\ t \leq \sqrt{N}}} \frac{1}{L(t)}. 
\end{split}\label{endinsight1}
\end{align}
 Since $L(t) > 0$ for each member $t$ in $S$, this and the third listed assertion 
 in Lemma \ref{fourpart} give us that 
\begin{align}
\begin{split}
 \sum_{\substack{t \in S \\ t \leq \sqrt{N}}} 
 \frac{1}{L(t)}
 & \leq \sum_{t \in S} \frac{1}{L(t)} \\ 
 & \leq \delta.
\end{split}\label{endinsight2}
 \end{align}
 Now, recall that $j \geq 1$ was arbitrary and that $N = N_{j}$. Also recall
 that $\theta = \frac{\varepsilon}{2}$ and that $\delta = \frac{\varepsilon}{4 K_{\theta}}$. 
 A combined application of the relations among 
 \eqref{endinsight1} and \eqref{endinsight2} gives us that 
\begin{align*}
 r_{\mathcal{C}}(N_{j}) 
 & = \frac{ \left| P(\mathcal{A}_{1}, \mathcal{A}_{2}) \cap [1, N_j] \right| }{N_j} \\ 
 & \leq \theta + K_{\theta} \delta \\ 
 & = \frac{\varepsilon}{2} + \frac{\varepsilon}{4} = \frac{3\varepsilon}{4} < \varepsilon. 
\end{align*}
 Since $(N_{j})_{j \geq 1}$ is a strictly increasing sequence of positive integers, 
 we find that $N_{j} \to \infty$, and we thus find that 
 $$ d_{\mathcal{C}} = \liminf_{M \to \infty} r_{\mathcal{C}}(M) 
 \leq \liminf_{j\to \infty} r_{\mathcal{C}}(N_{j}) \leq \frac{3\varepsilon}{4} < \varepsilon, $$ 
 as desired. 
\end{proof}

\begin{example}
 For  each $\varepsilon > 0$, our construction, as above,  gives us at least one coloring $\mathcal{C}$ of $\mathbb{N}$ such that 
 $d_{\mathcal{C}} < \varepsilon$. 
 There are many such $2$-colorings possible according to our construction, and this depends
 on how values of the form $M_{j}$, as in \eqref{firstMcondition}, 
 are chosen. Adopting the convention that $1$ is colored with red, our construction necessarily assigns colorings so that 
$$
\underbrace{1, 2, \ldots, M_1}_{\text{red}}, \ \underbrace{M_1 + 1, M_1 + 2, \ldots, M_1^2 + 1}_{\text{blue}}.
$$  For example, let $\varepsilon = \frac{1}{10}$. 
 In this case, the associated parameters are 
 $\theta = \frac{1}{20}$, 
 $K_{\theta} = 1060 + 48 \log 20$, and $\delta = \frac{1}{40 K_{\theta}}$. 
 By then writing 
 $ T = 84{\,}800 + 3840 \log 20$, 
 we define
 $$ M_{\min} = \min\left\{ m\in \mathbb{N} : m > 20, \sum_{p \leq m} \frac{1}{p} \geq T \right\}. $$
 The admissible choices for $M_1$ are thus $M_{\min}$, $M_{min} + 1$, etc. 
 Using classical bounds for sums of reciprocals of consecutive primes due to Rosser and Schoenfeld \cite[Theorem 5, Eq.\ (3.17)--(3.18)]{RosserSchoenfeld1962}, 
 it can be shown that 
$$ \exp(\exp(96{\,}303.3504)) 
 < M_{\min} \leq \exp(\exp(96{\,}303.3505)). $$
 So, in any case, the first $$\exp(\exp(96{\,}303.3504)) > 10^{10^{41{\,}823}}$$ integers would necessarily
 be colored red for the 
 $\varepsilon = \frac{1}{10}$ case. 
\end{example}

\section{Conclusion}
 Our solution to the Erd\H{o}s--S\'{a}rk\"{o}zy product-density problem indicated above and the techniques introduced and employed 
 above raise questions in regard to further areas of study related to the given problem and to our techniques. 
 For example, instead of coloring each positive \emph{integer} (red or blue) 
 and then considering lower limits as in 
 \eqref{definedcalC}, 
 we could instead color each \emph{prime} (red or blue)
 and then enumerate integers up to a given integer $N$
 that admit a monochromatic prime factorization. 
 This is decidedly inequivalent to the enumeration suggested in 
 \eqref{defineRCN}, and we encourage a full exploration
 of monochromatic prime factorizations in relation to the work of 
 Erd\H{o}s and S\'{a}rk\"{o}zy \cite{ErdosSarkozy1990,ErdosSarkozy1993}.

 For a given value $\alpha \in \big[\frac{1}{2}, 1\big]$, 
 how can we characterize 
 2-colorings $\mathcal{C}$ of $\mathbb{N}$ 
 such that the limit $\lim_{N \to \infty} r_{\mathcal{C}}(N)$ exists
 and, moreover, is equal to $\alpha$ (being consistent with the definition of $r_{\mathcal{C}}(N)$ in \eqref{defineRCN})?

 How can we characterize 2-colorings $\mathcal{C}$ of $\mathbb{N}$ such that 
 $\lim_{N \to \infty} r_{\mathcal{C}}(N)$ does not exist?

 \subsection*{AI disclosure}
 The following AI disclosure is based on AI disclosures used and written by the author previously (with numerous parts taken verbatim 
 from previous AI disclosures used and written by the author). 

 Our solution above to Guy's problem
 is, as indicated above, based on our extensive interactions with GPT-6 Astra Pro, and this is described in more 
 detail below. All 
 AI-generated suggestions were substantially revised, corrected, and independently verified by the author, who assumes full responsibility 
 for the mathematical content. The author's use of AI is fully in accordance with editorial standards on the responsible, ethical, and 
 transparent use of AI. 

 On 2026-10-18, the author inputted a prompt into GPT-6 Astra Pro containing an equivalent version of the problem from Guy's text 
 solved above (with an explicit attribution to Guy's text on unsolved problems in number theory), 
 concluding with an explicit request that GPT-6 Astra Pro solve the given problem. 
 After {\tt 13m 53s} of thinking, GPT-6 Astra Pro provided a response beginning with 
\begin{quote}
``{\bf I obtain a negative answer to the exact formulation you gave: no such absolute constant \(c>0\) exists.} The construction below produces, for every \(\varepsilon>0\), one fixed partition and an infinite sequence of cutoffs at which the proportion of distinct monochromatic products is at most \(\varepsilon\).''
\end{quote}

 After extensive clarifications, corrections, revisions, suggestions, verifications, calculations, etc., through extensive discussions between 
 the author and 
 GPT-6 Astra Pro, 
 and after extensive and careful reviews of such discussions on the part of the author, 
 this provided a basis for the above solution, written by the author, 
 to Guy's problem. 

 For example, it was through the author's enquiries, investigations, revisions, reviews, reflections, considerations, etc., 
 that the 
 ``non-probabilistic'' version of Markov's inequality 
 given above was employed to prove Lemma \ref{lemmaHbound}. Similar commentary 
 applies to the case analysis in the proof of Lemma \ref{lemmaHbound}, 
 the use of the Alzer--Berg inequality in \eqref{eqAlzerBerg}, 
 the expansion of $\sum_{n \leq M} (F(n) - L)^{2}$ involved in the proof of 
 Lemma \ref{lemmaHbound}, 
 the equivalent formulations of the 
 Erd\H os--S\'ark\"ozy problem discussed in Section \ref{secintro}, the evaluation in \eqref{evalan}, etc.

\bibliographystyle{plain}
\bibliography{seprefe}

@article {Hindman1979,
    AUTHOR = {Hindman, Neil},
     TITLE = {Partitions and sums and products of integers},
   JOURNAL = {Trans. Amer. Math. Soc.},
  FJOURNAL = {Transactions of the American Mathematical Society},
    VOLUME = {247},
      YEAR = {1979},
     PAGES = {227--245},
      ISSN = {0002-9947,1088-6850},
   MRCLASS = {10A45 (05A17 54A25)},
  MRNUMBER = {517693},
MRREVIEWER = {G.\ L.\ Alexanderson},
       DOI = {10.2307/1998782},
       URL = {https://doi.org/10.2307/1998782},
}

@article {RosserSchoenfeld1962,
    AUTHOR = {Rosser, J. Barkley and Schoenfeld, Lowell},
     TITLE = {Approximate formulas for some functions of prime numbers},
   JOURNAL = {Illinois J. Math.},
  FJOURNAL = {Illinois Journal of Mathematics},
    VOLUME = {6},
      YEAR = {1962},
     PAGES = {64--94},
      ISSN = {0019-2082},
   MRCLASS = {10.42},
  MRNUMBER = {137689},
       URL = {http://projecteuclid.org/euclid.ijm/1255631807},
}

@article {GreenSanders2016,
    AUTHOR = {Green, Ben and Sanders, Tom},
     TITLE = {Monochromatic sums and products},
   JOURNAL = {Discrete Anal.},
  FJOURNAL = {Discrete Analysis},
      YEAR = {2016},
     PAGES = {Paper No. 5, 48},
      ISSN = {2397-3129},
   MRCLASS = {11B30 (05D10)},
  MRNUMBER = {3533304},
MRREVIEWER = {Christian\ Elsholtz},
       DOI = {10.19086/da.613},
       URL = {https://doi.org/10.19086/da.613},
}

@incollection {Bergelson2005,
    AUTHOR = {Bergelson, Vitaly},
     TITLE = {Multiplicatively large sets and ergodic {R}amsey theory},
      NOTE = {Probability in mathematics},
   JOURNAL = {Israel J. Math.},
  FJOURNAL = {Israel Journal of Mathematics},
    VOLUME = {148},
      YEAR = {2005},
     PAGES = {23--40},
      ISSN = {0021-2172,1565-8511},
   MRCLASS = {37A45 (05D10)},
  MRNUMBER = {2191223},
MRREVIEWER = {Mathias\ Beiglb\"ock},
       DOI = {10.1007/BF02775431},
       URL = {https://doi.org/10.1007/BF02775431},
}

@incollection {ErdosSarkozy1990,
    AUTHOR = {Erd\H{o}s, P. and S\'{a}rk\"{o}zy, A.},
     TITLE = {On a conjecture of {R}oth and some related problems. {II}},
 BOOKTITLE = {Number theory ({B}anff, {AB}, 1988)},
     PAGES = {125--138},
 PUBLISHER = {de Gruyter, Berlin},
      YEAR = {1990},
      ISBN = {3-11-011723-1},
   MRCLASS = {11B05},
  MRNUMBER = {1106656},
MRREVIEWER = {H\'edi\ Daboussi},
}

@article {AlzerBerg2006,
    AUTHOR = {Alzer, Horst and Berg, Christian},
     TITLE = {Some classes of completely monotonic functions. {II}},
   JOURNAL = {Ramanujan J.},
  FJOURNAL = {Ramanujan Journal. An International Journal Devoted to the
              Areas of Mathematics Influenced by Ramanujan},
    VOLUME = {11},
      YEAR = {2006},
    NUMBER = {2},
     PAGES = {225--248},
      ISSN = {1382-4090,1572-9303},
   MRCLASS = {33B15 (11A41 26A48)},
  MRNUMBER = {2267677},
MRREVIEWER = {Yu.\ V.\ Trubnikov},
       DOI = {10.1007/s11139-006-6510-5},
       URL = {https://doi.org/10.1007/s11139-006-6510-5},
}

@article {ErdosSarkozy1993,
    AUTHOR = {Erd\H{o}s, P. and S\'{a}rk\"{o}zy, A.},
     TITLE = {Some solved and unsolved problems in combinatorial number
              theory. {II}},
   JOURNAL = {Colloq. Math.},
  FJOURNAL = {Colloquium Mathematicum},
    VOLUME = {65},
      YEAR = {1993},
    NUMBER = {2},
     PAGES = {201--211},
      ISSN = {0010-1354,1730-6302},
   MRCLASS = {11B13 (11B34)},
  MRNUMBER = {1240166},
MRREVIEWER = {A.\ J.\ Hildebrand},
       DOI = {10.4064/cm-65-2-201-211},
       URL = {https://doi.org/10.4064/cm-65-2-201-211},
}

@book {Rainville1960,
    AUTHOR = {Rainville, Earl D.},
     TITLE = {Special functions},
 PUBLISHER = {The Macmillan Company, New York},
      YEAR = {1960},
     PAGES = {xii+365},
   MRCLASS = {33.00},
  MRNUMBER = {107725},
MRREVIEWER = {R.\ W.\ Hamming},
}

@book {Guy1994,
    AUTHOR = {Guy, Richard K.},
     TITLE = {Unsolved problems in number theory},
    SERIES = {Problem Books in Mathematics},
   EDITION = {Second},
      NOTE = {Unsolved Problems in Intuitive Mathematics, I},
 PUBLISHER = {Springer-Verlag, New York},
      YEAR = {1994},
     PAGES = {xvi+285},
      ISBN = {0-387-94289-0},
   MRCLASS = {11-01 (00A07 11-02)},
  MRNUMBER = {1299330},
MRREVIEWER = {B.\ M. M. de Weger},
       DOI = {10.1007/978-1-4899-3585-4},
       URL = {https://doi.org/10.1007/978-1-4899-3585-4},
}

@article {Moreira2017,
    AUTHOR = {Moreira, Joel},
     TITLE = {Monochromatic sums and products in {$\Bbb N$}},
   JOURNAL = {Ann. of Math. (2)},
  FJOURNAL = {Annals of Mathematics. Second Series},
    VOLUME = {185},
      YEAR = {2017},
    NUMBER = {3},
     PAGES = {1069--1090},
      ISSN = {0003-486X,1939-8980},
   MRCLASS = {37A30 (05D10 11B75)},
  MRNUMBER = {3664819},
MRREVIEWER = {Siming\ Tu},
       DOI = {10.4007/annals.2017.185.3.10},
       URL = {https://doi.org/10.4007/annals.2017.185.3.10},
}

@article {CsikvariGyarmatiSarkozy2012,
    AUTHOR = {Csikv\'ari, P\'eter and Gyarmati, Katalin and S\'ark\"ozy,
              Andr\'as},
     TITLE = {Density and {R}amsey type results on algebraic equations with
              restricted solution sets},
   JOURNAL = {Combinatorica},
  FJOURNAL = {Combinatorica. An International Journal on Combinatorics and
              the Theory of Computing},
    VOLUME = {32},
      YEAR = {2012},
    NUMBER = {4},
     PAGES = {425--449},
      ISSN = {0209-9683,1439-6912},
   MRCLASS = {11B75 (05D10)},
  MRNUMBER = {2965285},
MRREVIEWER = {Nikolai\ Volodin},
       DOI = {10.1007/s00493-012-2697-9},
       URL = {https://doi.org/10.1007/s00493-012-2697-9},
}

@article {BeiglbockBergelsonHindmanStrauss2008,
    AUTHOR = {Beiglb\"ock, Mathias and Bergelson, Vitaly and Hindman, Neil
              and Strauss, Dona},
     TITLE = {Some new results in multiplicative and additive {R}amsey
              theory},
   JOURNAL = {Trans. Amer. Math. Soc.},
  FJOURNAL = {Transactions of the American Mathematical Society},
    VOLUME = {360},
      YEAR = {2008},
    NUMBER = {2},
     PAGES = {819--847},
      ISSN = {0002-9947,1088-6850},
   MRCLASS = {05D10},
  MRNUMBER = {2346473},
MRREVIEWER = {J.\ Spencer},
       DOI = {10.1090/S0002-9947-07-04370-X},
       URL = {https://doi.org/10.1090/S0002-9947-07-04370-X},
}

\end{document}